\documentclass[12pt]{amsart}
\usepackage{amssymb}
\usepackage[dvips]{graphics}
\ifx\mathrm\undefined\let\mathrm\rm\fi
\ifx\mathbf\undefined\let\mathbf\bf\fi
\ifx\mathfrak\undefined\let\mathfrak\frak\fi
\ifx\mathcal\undefined\let\mathcal\cal\fi
\ifx\mathbb\undefined\let\mathbb\Bbb\fi
\ifx\emph\undefined\let\emph\it\fi
\font\bb=msbm10 at9.98pt
\begin{document}
\def\semidirect{\hbox{$\;$\bb\char'156$\;$}}
\newcommand{\SL}{\mathrm{SL}}
\newcommand{\GL}{\mathrm{GL}}
\newcommand{\g}{{{\mathfrak g}\,}}
\newcommand{\bor}{{{\mathfrak b}}}
\newcommand{\n}{{{\mathfrak n}}}
\newcommand{\h}{{{\mathfrak h\,}}}
\newcommand{\Id}{{\operatorname{Id}}}
\newcommand{\Z}{{\mathbb Z}}
\newcommand{\ZZ}{{\mathbb Z_{>0}}}
\newcommand{\N}{{\mathbb N}}
\newcommand{\R}{{\mathbb R}}
\newcommand{\p}{{\mathbb P}}
\newcommand{\C}{{\mathbb C}}
\newcommand{\Q}{{\mathbb Q}}
\newcommand{\CC}{\mathcal{C}}
\newcommand{\HH}{\mathcal{H}}
\newcommand{\F}{\mathcal{F}}
\newcommand{\W}{\mathcal{W}}
\newcommand{\PP}{\mathcal{P}}
\newcommand{\Sym}{{\rm Sym}}
\newcommand{\Sing}{{\rm Sing}}
\newcommand{\Poly}{{\rm Poly}}
\newcommand{\Span}{{\rm Span}}
\newcommand{\Res}{{\rm Res}} 
\newcommand{\1}{{\bf 1}}
\newcommand{\kk}{{\bf k}}
\newcommand{\z}{{\bf z}}
\newcommand{\dontprint}[1]
{\relax}
\newtheorem%
{thm}{Theorem}
\newtheorem%
{prop}
{Proposition}
\newtheorem%
{lemma}
{Lemma}
\newtheorem%
{lemmadef}[thm]{Lemma-Definition}
\newtheorem%
{cor}
{Corollary}
\newtheorem%
{conj}
{Conjecture}
\newenvironment{definition}
{\noindent{\bf Definition\/}:}{\par\medskip}

\title  {Gaudin's model and the generating function\\ of the Wronski map}

\author[{}]
{I. Scherbak}

\maketitle     

\medskip
\centerline{\it School of Mathematical Sciences,
Tel Aviv University,}
\centerline{\it Ramat Aviv, Tel Aviv 69978, Israel}
\centerline{\it e-mail: \quad scherbak@post.tau.ac.il}

\bigskip

\pagestyle{myheadings}
\markboth{I. Scherbak}
{ Gaudin's model and the generating function of the Wronski map}  
\begin{abstract}
We consider  the  Gaudin model associated to a point $z\in\C^n$ with pairwise 
distinct coordinates and  to the subspace of singular vectors of a given weight 
in the tensor product of  irreducible finite-dimensional  $sl_2$-representations, 
\cite{G}. The Bethe equations of  this model  provide the critical point system
of a remarkable rational symmetric function.  Any critical orbit determines a common 
eigenvector of  the Gaudin hamiltonians called a  Bethe vector. 

In \cite{ReV}, it was shown that for generic $z$ the Bethe vectors  span the space 
of singular vectors, i.e. that  the number of  critical orbits is bounded from below
by the dimension of the space of singular vectors. The upper bound by the same
number is one of the main results of  \cite{SV}.

In the present paper we get  this upper bound in another, ``less technical'',  way.
The crucial observation is that the symmetric function  defining the Bethe equations
can be interpreted as the generating  function of the  map sending a pair of complex 
polynomials into their Wronski determinant: the  critical orbits
determine the preimage of a given polynomial under this map. 
Within the  framework of the Schubert calculus, the number of  critical orbits  can be 
estimated by the  intersection  number  of  special Schubert classes.  
Relations to the  $sl_2$ representation theory  (\cite{F}) imply that this number 
is  the dimension of the space of singular vectors.

We prove also that the spectrum of the Gaudin hamiltonians is simple for generic $z$. 
\end{abstract}

\section{Introduction}\label{S1}
The Gaudin model of statistical mechanics is a completely integrable quantum spin  
chain associated to the Lie algebra  $sl_2=sl_2(\C)$, \cite{G}. 
Denote  $L_{\lambda_j}$ the  irreducible  $sl_2$-module  with  highest weight 
$\lambda_j\in\C$. The space of states of the model is the tensor product
$L=L_{\lambda_1}\otimes \dots \otimes  L_{\lambda_n}$. 
Associate with any  $L_{\lambda_j}$ a complex number $z_j$, and  assume 
$z_1,\dots,z_n$ to be pairwise distinct numbers, $z=(z_1,\dots,z_n)$.
The Gaudin hamiltonians $\HH_1(z), \dots, \HH_n(z)$ are  mutually commuting
linear operators on $L$ which are defined as follows, 
$$
\HH_j(z)=\sum_{i\neq j}\frac{\Omega_{ij}}{z_j-z_i}\,, \quad  j=1,\dots, n,
$$
where $\Omega_{ij}$ is the  operator which acts  as the Casimir element on $i$-th and $j$-th 
factors of  $L$ and as the identity on all others. 

\medskip
One of the main problems  in the Gaudin model is to find common eigenvectors
and eigenvalues of the Gaudin hamiltonians, and the algebraic Bethe Ansatz 
is  one of  the most effective  methods  for solving this problem.
The idea of this method is to find some function 
with values in the space of states, and to determine a certain special value 
of its argument  in such a way  that the value of this function is an eigenvector. 
The  equations which determine these  special values of the argument are called 
{\it the Bethe equations}, and the common eigenvector  corresponding 
to  a solution of the Bethe equations,  is called {\it the Bethe vector}, \cite{FaT}. 

\medskip
It is enough to diagonalize the Gaudin hamiltonians 
in the subspace of singular vectors of a given weight.
The Bethe equations associated to the space $\Sing_k(L)$ of singular vectors of the weight
$\lambda_1+\dots +\lambda_n -2k$, where $k$ is a positive integer,  have the form
$$
\sum_{l=1}^n\frac{\lambda_l}{t_i-z_l}\,=\,\sum_{j\neq i}\frac{2}{t_i-t_j}\,,
\quad  i=1,\dots, k\,.
$$
This system is symmetric with respect to permutations of  unknowns $t_1,\dots,
t_k$,  and any  (orbit of)  solution $t^0=(t_1^0,\dots, t_k^0)$ to this system defines 
an eigenvector $v(t^0)\in\Sing_k(L)$ of  the Gaudin hamiltonians $\{\HH_j(z)\}$, \cite{G}.

\medskip
A conjecture related to Bethe Ansatz says that the  Bethe vectors  give a basis 
of the space of states.  For the $sl_2$ Gaudin model,
the conjecture was proved in \cite{ReV} in the case of generic $z$  and  generic 
(non-resonant) weights $\lambda_1, \dots,\lambda_n$, 
In the case of {\it integral dominant} weights, which is a non-generic one, 
results of Sec.~ 9 of  \cite{ReV}  imply only that for generic $z$   the number of 
Bethe vectors is {\it at least} the dimension of the space of singular vectors, 
see also Theorem~8 in \cite{SV}.  The proof of the Bethe Ansatz conjecture 
for the $sl_2$ Gaudin model was completed in our work with A.~Varchenko \cite{SV},
where the bound from above for the number of the Bethe vectors by the dimension
of the space of singular vectors was obtained.  

\medskip
In \cite{SV}, we study the function 
$$
\Phi(t)=\prod_{i=1}^k\prod_{l=1}^n (t_i-z_l)^{-\lambda_l}
\prod_{1\leq i<j\leq k}(t_i-t_j)^2\,. 
$$
If all $\lambda_l$ are positive integers, then this function is  rational,
and the Bethe equations are exactly the equations on the critical
points of this function with non-zero critical values,
$$
\frac 1{\Phi(t)}\,\cdot\,\frac {\partial \Phi}{\partial t_i}(t)=0\,\quad i=1,\dots, k.
$$
The upper bound for the number of critical orbits of the function $\Phi(t)$
with non-zero critical values by the dimension of $\Sing_k(L)$ 
is one of the main results of  the paper. The proof  is difficult, 
see Theorems 9--11 of  \cite{SV}.  As we point out in  Sec.~1.4 of \cite{SV}, 
the critical orbits are labeled  by  certain two-dimensional planes in the linear space 
of complex polynomials.  This observation suggests to apply the Schubert calculus 
to this problem. In the present paper we realize this approach.

\medskip
The Wronski determinant of  two polynomials  in one variable defines 
a map from the Grassmannian of two-dimensional  planes of the linear
space of complex polynomials  to the space of monic polynomials called
{\it the Wronski map}.

\medskip
The function $\Phi(t)$ turns out to be the generating function of the Wronski
map:  the critical orbits  label planes in the preimage under this map 
of the polynomial
$$
W(x)=(x-z_1)^{\lambda_1}\dots (x-z_n)^{\lambda_n}\,.
$$
In fact, this is a reformulation of a classical result going back to Heine and Stieltjes, 
Ch.~6.8 of  \cite{Sz}. 
 
 \medskip
To calculate the cardinality of the preimage of the Wronski map 
is a problem of enumerative algebraic geometry, and an upper bound can be 
easily obtained in terms of the intersection  number of special Schubert 
classes. A well-known relation between  representation 
theory and  the Schubert calculus (\cite{F}) implies that the obtained upper 
bound coincides with the dimension  of  $\Sing_k(L)$.

\medskip
 In Sec.~\ref{S2} we collect known facts
related to the Gaudin model.  In Sec.~\ref{S3} we show that the function
defining the Bethe equations is the generating function of the Wronski map,
and in Sec.~\ref{S4} we obtain an upper bound for the number of critical
orbits of the generating function in terms of the intersection number of special
Schubert classes.

\medskip
Another  conjecture related to the Gaudin model says that for generic $z$
the Gaudin hamiltonians have a simple spectrum. 
In \cite{ReV} it was proved that Bethe vectors are differ by their eigenvalues
for  generic  $\lambda_1, \dots,\lambda_n$  and  real  $z$ of the form $z_j=s^j$, 
where $s>>1$.  If it is known that the Bethe vectors form a basis, then the 
simplicity of the spectrum follows.

\medskip
In Sec.~\ref{S5} we deduce the simplicity of  the spectrum of the Gaudin hamiltonians  
for  generic $z$  and integral dominant  $sl_2$ weights from the relation to Fuchsian 
differential equations described in \cite{SV}.

\medskip
The argumets presented here work and give the similar results for the $sl_p$ Gaudin
model associated with the tensor product of symmetric powers of the standard
$sl_p$-representation, see Sec.~5 of \cite{S1}. The link between the Bethe vectors and
the Wronski maps seems to be useful  for stufy rational maps;  see \cite{S2}
where the case of rational functions is treated.

\medskip
The author is thankful to A.~Eremenko, E.~Frenkel, F.~Sottile, A.~Varchenko, V.~Zakalyukin 
for useful  discussions, to G.~Fainshtein and D.~Karzovnik for stimulating comments.

\section{Gaudin's model}\label{S2}
\subsection{The Gaudin hamiltonians}\label{s21}
Let   $m_1,\dots, m_n$ be nonnegative integers,  $M=(m_1,\dots, m_n)$.
Denote  $L_{m_j}$ the  irreducible  $sl_2$-module  with  highest weight $m_j$.
{\it The space of states} of the Gaudin model is the tensor product 
\begin{equation}
\label{LM}
L^{\otimes M}=L_{m_1}\otimes \dots \otimes  L_{m_n}\,.
\end{equation}

\medskip
Associate with any  $L_{m_j}$ a complex number $z^0_j$, and  assume 
$z^0_1,\dots,z^0_n$ to be pairwise distinct numbers, $z^0=(z^0_1,\dots,z^0_n)$. 
Let  $e, f, h$  be the standard generators of $sl_2$,
$$
[e,f]=h,\quad  [h,e]=2e,\quad [h,f]=-2f
$$
and                                  
$$
\Omega=e\otimes f+f\otimes e+\frac12 h\otimes h \in sl_2\otimes sl_2
$$
 the Casimir element.
For $1\leq i<j\leq n$, denote $\Omega_{ij}:L^{\otimes M}\to L^{\otimes M}$  
the operator which acts  as $\Omega$ on $i$-th and $j$-th factors of
$L^{\otimes M}$ and as the identity on all others. 
{\it The  hamiltonians} of the  Gaudin model are  defined as follows,
\begin{equation}
\label{H}
\HH_i(z^0)=\sum_{j\neq i}\frac{\Omega_{ij}}{z^0_i-z^0_j}\,, \quad  i=1,\dots, n.
\end{equation}

\subsection{Subspaces of singular vectors}\label{s23}
Write $ |M|=m_1+\dots +m_n$. For  a nonnegative integer $k$ 
such that $|M|-2k\geq 0$, define $\Sing_k\,$, 
the subspace of  singular vectors of weight  $|M|-2k$ in $L^{\otimes M}$,
\begin{equation}
\label{Sing}
\Sing_k\,=\, \Sing_k (L^{\otimes M})=
\left\{\, w\in  L^{\otimes M}\ \vert\  ew=0\,,\ hw=(|M|-2k)w\,\right\}.
\end{equation}

\begin{thm}  {\rm (\cite{SV}, Theorem~5 )}\label{tI} We have
$$
\dim\Sing_k(L^{\otimes M})\,=\,\sum_{q=0}^{n} (-1)^q
\sum_{1\leq i_1<\dots <i_q\leq n}
{k+n-2-m_{i_1}-\dots -m_{i_q}-q\choose n-2}\,,
$$
where we assume  ${a\choose b}=0$ for $a<b$.  \hfill $\triangleleft$
\end{thm}

\medskip
The subspace  $\Sing_k $ is an invariant subspace of the Gaudin hamiltonians
 for  any $0\leq k\leq |M|/2$,  and the singular vectors generate the whole of 
 $L^{\otimes M}$. Therefore it is enough to diagonalize Gaudin hamiltonians 
 in a given $\Sing_k$. 

\subsection{Bethe equations associated to  $\Sing_k $ and $z^0$}\label{s24} 
The  Bethe equations for the Gaudin model associated
to  $\Sing_k $ and $z^0$ have the form
\begin{equation}\label{cps}
\sum_{l=1}^n\frac{m_l}{t_i-z^0_l}\,=\,\sum_{j\neq i}\frac{2}{t_i-t_j}\,,
\quad  i=1,\dots, k\,.
\end{equation}
Any  solution $t^0=(t_1^0,\dots, t_k^0)$ to this system define an eigenvector
$v(t^0)\in\Sing_k$ of  the Gaudin hamiltonians, $\HH_j(z^0)v(t^0)=\mu_jv(t^0)$, 
with eigenvalues
\begin{equation}\label{cvs}
\mu_j=\mu_j(t^0)=\sum_{l\neq j}\frac{m_lm_j}{2(z^0_j-z^0_l)}\,-\,
\sum_{i=1}^k\frac{m_j}{z^0_j-t^0_i}\,,
\quad  j=1,\dots, n\,.
\end{equation}

\subsection{The master function of the Gaudin model}\label{s25}
Consider   the following function in  variables
$t=(t_1,\dots,t_k)$ and $z=(z_1,\dots,z_n)$,
\begin{equation}\label{Psi}
\Psi(t,z)=\prod_{1\leq i<j\leq n}(z_i-z_j)^{m_im_j/2}
\prod_{i=1}^k\prod_{l=1}^n (t_i-z_l)^{-m_l}
\prod_{1\leq i<j\leq k}(t_i-t_j)^2\,, 
\end{equation}
defined on 
\begin{eqnarray}\label{CC}
\CC&=&\CC(t;z)= \   \left\{\, t\in\C^k\,,  z\in\C^n\,  \vert\ \right.\\
\ &\ &\left.t_i\neq z_j,\, t_i\neq t_l,\, z_j\neq z_q,\, 1\leq i\neq l\leq k,\, 
1\leq j\neq q\leq n\, \right\}\nonumber\,.
\end{eqnarray}
Write 
\begin{eqnarray*}
S(t,z)& = &\ln \Psi(t,z)=\\
\ & \ & \sum_{1\leq i<j\leq n}\frac{m_im_j}2 \ln (z_i-z_j)
-\sum_{i=1}^k\sum_{l=1}^n m_l\ln (t_i-z_l) +\sum_{1\leq i<j\leq k}2\ln (t_i-t_j).
\end{eqnarray*}
The both functions  $S(t,z)$ and $\Psi(t,z)$ have clearly the same critical
set in $\CC(t;z)$. Denote 
$$
 Z=\left\{\,  z\in\C^n\,  \vert\, z_j\neq z_q,\  1\leq j\neq q\leq n\, \right\}.
$$
Let $\pi\,:\, \CC\to Z,\quad \pi(t,z)= z,$ be the natural projection.
For fixed $z^0\in Z$, 
the equations (\ref{cps})  form a  critical point system of the function $S(t,z^0)$
considered as a function in $t$  on $\pi^{-1}(z^0)$,
$$
\frac{\partial S}{\partial t_i}(t,z^0)=0\,, \
i=1,\dots, k\,,
$$
and the equalities  (\ref{cvs})  can be re-written in the form
$$
\mu_j=\frac{\partial S}{\partial z_j}(t^0,z^0)\,, \ j=1,\dots,n,
$$
where $t^0=t^0(z^0)$ is a solution to the Bethe equations. 
In terms of the function $\Psi(t,z)$, the Bethe equations (\ref{cps}) are
$$
\frac1{\Psi(t,z^0)}\,\cdot\,\frac{\partial\Psi}{\partial t_i}(t,z^0)=0\,, \
i=1,\dots, k\,,
$$
and
$$
\mu_j=\frac1{\Psi(t^0,z^0)}\,\cdot\,\frac{\partial\Psi}{\partial z_j}(t^0,z^0)\,,
\ j=1,\dots,n.
$$

\subsection{The Shapovalov form}\label{s26}
Let $m$ is a nonnegative integer.
For the $sl_2$-module $L_{m}$ with highest weight $m$, fix  the highest weight singular  vector
$$
v_m\in L_{m}\,,\quad  hv_m=mv_m,\quad ev_m=0.
$$
Denote $B_m$   the unique bilinear symmetric form  on $L_{m}$ such that
$$
B_m(v_m,v_m)=1,\quad B_m(hx,y)=B_m(x,hy)\,,\quad  B_m(ex,y)=B_m(x,fy)
$$ 
for all $x,y\in L_{m}$. The   vectors 
$v_m, fv_m,\dots, f^{m}v_m$ are orthogonal with respect to $B_m$ and
form a basis of  $L_{m}$.

\medskip
The  bilinear symmetric form  on  $L^{\otimes M}$  given by  
 \begin{equation}
 \label{B}
B=B_{m_1}\otimes\dots\otimes B_{m_n}
\end{equation}
is  called {\it the Shapovalov form}.

\medskip
Let $j_1,\dots,j_n$ be  integers such that $0\leq j_i\leq m_i$
for any $1\leq i\leq n$. Write $J=(j_1,\dots,j_n)$ and $|J|=j_1 +\dots +j_n\,$. Denote
\begin{equation}\label{fJ}
f^Jv_M=f^{j_1}v_{m_1}\otimes \dots \otimes f^{j_n}v_{m_n}\,.
\end{equation}
The vectors $\left\{f^Jv_M\right\}$ are orthogonal with respect to the Shapovalov  
form  $B$ and provide a basis of  the space $L^{\otimes M}$.  We have 
$$
h(f^Jv_M)\,=\,(|M|-2|J|)f^Jv_M\,,\quad e(f^Jv_M)=0\,,
$$
i.e. the vector $f^Jv_M$ is {\it a singular vector of weight} $|M|-2|J|$.
The space $\Sing_k$ is generated by the vectors $f^Jv_M$  with $|J|=k$.

\subsection{Bethe vectors}\label{s27}
For $J=(j_1,\dots,j_n)$ with integer coordinates such that
$0\leq j_i\leq m_i$ and  $|J|=k$ and for $(t,z)$ in $\CC$ given by (\ref{CC})  set
$$
A_J(t,z)=\frac1{j_1!\dots j_n!}\Sym_t\left[\, \prod_{l=1}^n
\prod_{i=1}^{j_l}\frac1{t_{j_1+\dots+j_{l-1}+i}-z_l}\, \right]\,,
$$
where  
$$
\Sym_tF(t)=\sum_{\sigma\in S^k}F\left(t_{\sigma(1)},\dots, t_{\sigma(k)}\right)
$$
is the sum  over all  permutations  of $t_1,\dots, t_k$.

\begin{thm}\label{t2}{\rm (\cite{ReV})}\ \   (i)\ \
If $t^{(i)}$ is a nondegenerate critical point of  the function $S(t,z^0)$, then the vector
\begin{equation}\label{Bv}
v(t^{(i)},z^0)=\sum_{J\, :\, |J|=k} A_J(t^{(i)},z^0)f^Jv_M\,
\end{equation}
is an eigenvector of the operators $\HH_1(z^0),\dots, \HH_n(z^0)$.

\noindent
(ii)\ \
For generic $z^0$, the eigenvectors $v(t^{(i)},z^0)$ corresponding to
all critical point of  the function $S(t,z^0)$ generate $\Sing_k\,$. 
 \hfill $\triangleleft$
\end{thm}

The words ``generic $z^0$'' mean that $z^0$ does not belong to a suitable
proper algebraic set of $\C^n$.

\medskip
The set of  critical points of the function  $S(t,z^0)$ is invariant
with respect to the permutations of $t_1,\dots,t_k$, and
critical points belonging to the same orbit clearly define the same vector.
Theorem~\ref{t2} gives a lower bound for the number of critical orbits of
the function  $S(t,z^0)$ by the dimension of  $\Sing_k\,$. 
In \cite{SV} the upper bound by the same number was obtained.

\begin{thm}\label{t3}  {\rm \cite{SV}}  For fixed generic $z^0$,
all critical points of the  function   $S(t,z^0)$  are nondegenerate,
and the number of orbits of critical points is at most the dimension of  
$\Sing_k\,$.
 \hfill $\triangleleft$
\end{thm}

\begin{cor}{\rm \cite{SV}}
For generic $z^0$, the Bethe vectors of the $sl_2$ Gaudin model form 
a basis in the space of singular vectors of  a given weight. 
 \hfill $\triangleleft$
\end{cor}

\medskip
The statement that  all critical points of the  function   $S(t,z^0)$  
are nondegenerate is easy, see  Theorem~6  in  \cite{SV}. 
The difficult part of  \cite{SV} was to  estimate from above  the number of  critical
orbits  (see Theorems 9--11 in  \cite{SV}).  As we pointed out in  Sec.~1.4 of \cite{SV}, 
the orbits of  critical points  of  the function  $S(t,z^0)$ are labeled  by  certain 
two-dimensional planes in the linear space of complex polynomials. 
This observation suggests to apply the Schubert calculus in order to prove the
following

\medskip\noindent
{\bf  Proposition.}\ \ 
{\it The number of orbits of critical points of the function   $S(t,z^0)$
 is at most the dimension of $\Sing_k\,$.}

\medskip
As we will see below, within the framework of the Schubert calculus, 
the Proposition  is a direct and immediate corollary of easy Theorem 
\ref{tWr} of Sec.~\ref{s42}. 

\section{The generating function of the Wronski map}\label{S3}
\subsection{The Wronski map}\label{s31}
Let  $\Poly_d$ be  the vector space of  complex polynomials of degree 
at most $d$ in one variable. Denote  $G_2(\Poly_d)$ the Grassmannian of  
two-dimensional  planes in  $\Poly_d$. The complex dimension of    
$G_2(\Poly_d)$ is $2d-2$.

 \medskip
 For any $V\in G_2(\Poly_d)$ define  {\it the degree } 
of $V$ as  the maximal degree  of its polynomials and the {\it order}  of $V$ as  
the minimal   degree of its non-zero polynomials.
Let $V\in G_2(\Poly_d)$ be a plane of order $a$ and of degree $b$. 
Clearly $0\leq a<b\leq d$. Choose in $V$ two  monic polynomials, $F(x)$ and $G(x)$,
 of  degrees   $a$ and $b$ respectively.  They form  a basis of $V$.   
{\it The Wronskian of} $V$  is defined as the monic polynomial
$$
W_V(x)=\frac{F'(x)G(x)-F(x)G'(x)}{a-b}\,.
$$

The following lemma is evident.
\begin{lemma}\label{l1}
(i)\ \ The degree of $W_V(x)$ is  $a+b-1\leq 2d-2$.

\noindent
 (ii)\ \ The   polynomial  $W_V(x)$  does not depend on the  choice of a monic basis.

  \noindent
 (iii)\ \ All polynomials  of   degree $a$  in $V$ are proportional.
\hfill $\triangleleft$
\end{lemma}

 Thus  the mapping sending $V\in G_2(\Poly_d)$ to $W_V(x)$ is  a well-defined map from   
 the Grassmannian $G_2(\Poly_d)$ to $\C\p^{2d-2}$. We call it  {\it the Wronski map}. 
 This is a mapping between 
 smooth complex algebraic varieties of the same dimension,  and  hence the preimage
 of  any polynomial consists of  a finite number of   planes. 
On Wronski maps see \cite{EGa}.
 
 \subsection{Planes with a given Wronskian}\label{s32}
\begin{lemma}\label{l2}
 Any element of  $ G_2(\Poly_d)$ with a given Wronskian is uniquely determined 
 by any  of its  polynomial. 
\end{lemma}

\noindent
{\bf Proof:}\  \  Let $W(x)$ be the Wronskian of a plane $V\in G_2(\Poly_d)$ 
and  $f(x)\in V$. Take any  polynomial  $g(x)\in V$ linearly independent with $f(x)$. 
The plane $V$  is the solution space of the following second order   linear
differential equation  with respect to  unknown function   $u(x)$,
$$
\left|\begin{array}{ccc}u(x)\ \ &f(x)\   \ &g(x)\\
                          u'(x)\ \ &f'(x)\ \  &g'(x)\\
                          u''(x)\  \ &f''(x)\ \ &g''(x)
			   \end{array}\right|=0.
$$
The Wronskian of the polynomials $f(x)$ and $g(x)$ is proportional 
to $W(x)$, therefore this equation can be re-written in the form 
$$
W(x)u''(x)-W'(x)u'(x)+h(x)u(x)=0\,,
$$
where 
$$
h(x)\, =\, \frac{-W(x)f''(x)+W'(x)f'(x)}{f(x)}\,,
$$ 
as $f(x)$ is clearly  a solution to this equation.
\hfill $\triangleleft$

\medskip
We  call a plane $V\in G_2(\Poly_d)$  {\it generic} if for any $x_0\in\C$  
there is a polynomial  $P(x)\in V$ such that  $P(x_0)\neq 0$. In a generic plane, 
the polynomials of  any basis do not have common roots, and  almost all polynomials  
of the bigger degree do not have multiple roots. 
 
\medskip
The following lemma is evident.
\begin{lemma}\label{l3}  Let $V\in G_2(\Poly_d)$  be  a generic plane.

 \noindent
(i)\ \  If  $P(x)=(x-x_1)\dots(x-x_l)\in V$ is a polynomial  without
multiple roots, then $W_V(x_i)\neq 0$ for all $1\leq i\leq l$. 

 \noindent
(ii)\ \  If   $x_0$  is   a root of multiplicity  $\mu>1$ of  a
polynomial  $Q(x)\in V$, then  $x_0$ is a root of  $W_V(x)$
of multiplicity  $\mu-1$.  \hfill $\triangleleft$
\end{lemma}

A generic plane $V$ is {\it nondegenerate} if  the polynomials  of  the smaller degree 
 in $V$ do not have multiple roots.

\begin{lemma}\label{l4} Let $V\in G_2(\Poly_d)$ be a nondegenerate plane. 
If  the Wronskian $W_V(x)$ has the form 
$$
W_V(x)=x^{m}\tilde W(x)\,,\ \  \tilde W(0)\neq 0\,,
$$
then there exists a polynomial   $F_0(x)\in V$ of the form  
$$
F_0(x)=x^{m+1}\tilde F(x)\,,\ \  \tilde F(0)\neq 0\,.
$$
\end{lemma}

\noindent
{\bf Proof:}\ \  Let $G(x)\in V$ be a polynomial of the smaller degree.
We have $G(0)\neq 0$, due  to Lemma \ref{l3}. Let
 $F(x)\in V$  be a polynomial of the bigger degree. The polynomial
$$
F_0(x)\,=\,F(x)-\frac{F(0)}{G(0)}\,G(x)\,\in V
$$
satisfies  $F_0(0)=0$  and therefore  has the form
$$
F_0(x)=x^l  \tilde F(x)\,, \ \ \tilde F(0)\neq 0\,,
$$
for some integer  $l\geq 1$.  The polynomials $G$ and $F_0$ form a basis
of $V$, therefore  the polynomial $F'_0(x)G(x)-F_0(x)G'(x)\,$ is 
proportional to $W_V(x)$.  The smallest degree term in
this polynomial is  $l\cdot a_l\cdot G(0)x^{l-1}\,$,
where  $a_l$ is the coefficient of $x^l$ in $F_0(x)$.  Therefore $l=m+1$.
 \hfill $\triangleleft$

\subsection{Nondegenerate planes in $\Poly_d$ with a given Wronskian}\label{s321}
Recall that {\it the resultant} ${\rm Res}(P,Q)$ of polynomials  
$P(x)$ and  $Q(x)$  is an irreducible  integral polynomial in the coefficients 
of $P(x)$ and $Q(x)$ which vanishes whenever  $P(x)$ and $Q(x)$ 
have a common root,   and {\it the discriminant}  $\Delta(P)$  of $P(x)$  
is an irreducible integer  polynomial in the  coefficients  of $P(x)$ 
which vanishes whenever  $P(x)$ has a multiple root.

\medskip
Our aim is to describe the nondegenerate planes of a given order in the preimage 
of the polynomial   
\begin{equation}
\label{W}
W(x)=(x-z_1)^{m_1}\dots (x-z_n)^{m_n}\,.
\end{equation}
under the Wronski map.

\medskip
If plane $V$ of order $k$ and of degree $d$ has the Wronskian $W(x)$,
then $d=|M|+1-k>k$.
Let $F(x)$ be an unknown  polynomial of degree $1\leq k\leq |M|/2$.
Consider the following function, which appeared in unpublished notes
of V.~Zakalyukin related to the Wronski map,
$$
 \Phi\,=\,\Phi (F; W)=\,\frac{\Delta(F)}{\Res(W,F)}\,.
$$
If the polynomial $F(x)$ belongs to a plane in $\Poly_{|M|+1-k}$
with the Wronskian $W(x)$,  then $W(x)=F'(x)Q(x)-F(x)Q'(x)$ for some polynomial 
$Q(x)$.  Differentiating gives $W'(x)=F''(x)Q(x)-F(x)Q''(x)\,$, and we have 
$$
\frac{W'}{W}\,=\,\frac{F''Q-FQ''}{F'Q-FQ'}\,.
$$
If the plane spanned by $F(x)$ and $G(x)$ is  non-degenerate, then
polynomials $F$ and $Q$ do not have common roots,  polynomial $F$ does not
have multiple roots, and polynomials $W$ and $F$ do not have common roots,
by Lemmas~\ref{l3} and \ref{l4}.
Therefore at each root $t_i$  of $F$ we get
\begin{equation}\label{wf}
\frac{W'(t_i)}{W(t_i)}\,=\,\frac{F''(t_i)}{F'(t_i)}\,,\quad  i=1,\dots,k.
\end{equation}

\subsection{Theorem of  Heine--Stieltjes}\label{s33}
In order to re-write the function $\Phi$ as a function in unknown roots of  $F(x)$,
recall that if  $F(x)=(x-t_1)\dots (x-t_k)$  and if $W(x)$ is as in (\ref{W}), then
$$
\Delta(F)=\prod_{1\leq i<j\leq k}(t_i-t_j)^2\,,\quad
\Res(W,F)=\prod_{i=1}^k\prod_{j=1}^n(t_i-z_j)^{m_j}\,.
$$ 
We have
\begin{equation}\label{Phi}
\Phi=\Phi(t;z,M)=\prod_{i=1}^k\prod_{l=1}^n (t_i-z_l)^{-m_l}
\prod_{1\leq i<j\leq k}(t_i-t_j)^2\,. 
\end{equation}

The critical points of this function  were studied
by Heine and Stieltjes  in connection with  second
order linear differential equations having polynomial coefficients and a polynomial
solution of a prescribed degree. The result  of Heine and Stieltjes can be formulated 
as follows.

\begin{thm}\label{tW} {\rm (Heine--Stieltjes, cf. \cite{Sz},  Ch.~6.8)}\quad Let 
$t^0$  be a critical point with non-zero critical value of  the function 
$\Phi(t)$ given by {\rm (\ref{Phi})}.
Then $F(x)=(x-t^0_1)\dots(x-t^0_k)$ is  a polynomial of the smaller degree
in a nondegenerate plane with  the Wronskian $W(x)$ given by  {\rm (\ref{W})}.

Conversely,  if $F(x)=(x-t^0_1)\dots(x-t^0_k)$ is a  polynomial  of the smaller
degree in  a non-degenerate plane $V$ with the Wronskian $W_V(x)=W(x)$, then 
 $t^0=(t^0_1,\dots, t^0_k)$ is a critical point with non-zero critical value of  
 the function $\Phi(t)$.
\hfill $\triangleleft$   
\end{thm}
The function $\Phi(t)$ is symmetric with respect to permutations of $t_1,\dots, t_k$,
critical points belonging to one orbit define the same polynomial $F(x)$ and hence,
according to Lemma~\ref{l2}, the same plane.

\subsection{Bethe vectors and  nondegenerate planes}\label{s34}
The function $\Phi(t)$ given by (\ref{Phi})  and  
the function $\Psi(t,z)$ given by (\ref{Psi} are differ in
a factor depending only on $z$ and $M$,
$$
\Psi(t,z)=\prod_{1\leq i<j\leq n}(z_i-z_j)^{m_im_j/2}\cdot\,\Phi(t;z,M)\,.
$$
An easy calculation shows that  (\ref{wf}) is exactly the system defining the critical 
points with  non-zero critical values of the function $\Phi$, considered as a function in $t$. 
The theorem of  Heine--Stieltjes implies the following result.

\begin{cor}\label{c2}
There is a one-to-one correspondence 
between the critical orbits with non-zero critical value of  the  function   
{\rm (\ref{Phi})}  
and the nondegenerate planes  of  order $k$  and of degree $|M|+1-k$  having   
Wronskian {\rm (\ref{W})}.   \hfill $\triangleleft$
\end{cor}

\begin{cor}
For fixed $z\in Z$,  $M=(m_1,\dots, m_n)$ and $k$, there is a one-to-one correspondence 
between the Bethe vectors of the $sl_2$ Gaudin model associated to $z$ and $\Sing_k$ 
given by (\ref{Sing}) and nondegenerate planes of order $k$  with the  Wronskian  (\ref{W}). 
 \hfill $\triangleleft$
\end{cor}
 
 \begin{cor}
The Bethe equations associated to  $z$ and  {\rm (\ref{Sing})}  coincide 
with the critical point system of the generating function associated to order  $k$  and  
the Wronskian  {\rm (\ref{W})}.  \hfill $\triangleleft$ 
\end{cor}

The function $\Phi$ is called {\it the generating function} 
of the Wronski map: for any given monic polynomial $W(x)$ and any given order 
$k\leq \deg W/2$,  the  critical orbits of the  function  $\Phi$ determine  
the non-degenerate planes of order $k$ in the preimage of $W(x)$.  
We will show below  that for generic $z$ the critical orbits
determine {\it all} planes of order $k$ in the preimage, see Corollary~\ref{last} 
in Sec.~\ref{s42}. 

\section{The preimage of a given Wronskian}\label{S4}
The number of  nondegenerate planes  with a given  Wronskian 
can be estimated from above by the intersection number of Schubert classes.

\subsection{Schubert calculus {\rm (\cite{GrH}, \cite{F})}}\label{s41}
Let   $G_2(d+1)= G_2(\C^{d+1})$ be the Grassmannian variety of  
two-dimensional subspaces $V\subset \C^{d+1}$.   
A chosen basis  $e_1,\dots,e_{d+1}\,$ of $ \C^{d+1}$
defines the flag of linear subspaces
$$
E_{\bullet}\,:\quad E_1\, \subset \, E_2\, \subset\, \dots\,  
\subset \, E_d \subset \, E_{d+1}= \C^{d+1}\,,
$$
where $E_i=\Span\{e_1,\dots,e_i\}\,$, $\dim E_i=i$.   
For any integers  $a_1$ and $a_2$
such that   $0\leq a_2\leq a_1\leq d-1\,$, the {\it Schubert variety} 
$\Omega_{a_1,a_2}(E_{\bullet})\subset  G_2({d+1})$
is defined as follows,
$$
\Omega_{a_1,a_2}= \Omega_{a_1,a_2}(E_{\bullet})=\left\{\, V\in  G_2({d+1})\, 
\vert\, \dim \left(V\cap E_{d-a_1}\right)\geq 1\,,
\   \dim \left(V\cap E_{d+1-a_2}\right)\geq 2\, \right\}\,.
$$
The variety $\Omega_{a_1,a_2}=\Omega_{a_1,a_2}(E_{\bullet})$
is an irreducible closed subvariety of $ G_2({d+1})$  of the complex
codimension $a_1+a_2$.

\medskip
The homology classes $[\Omega_{a_1,a_2}]$  of Schubert varieties
 $\Omega_{a_1,a_2}$ are independent of the choice of  flag, and form a basis
for the integral homology of  $G_2(d+1)$.
Denote  $\sigma_{a_1,a_2}$  the cohomology class in 
 $H^{2(a_1+a_2)}(G_2({d+1}))$ whose cap product with the
fundamental class of   $G_2(d+1)$ is the homology class
 $[\Omega_{a_1,a_2}]$.
The classes  $\sigma_{a_1,a_2}$ are called {\it Schubert classes}.
They give a basis over $\Z$ for the cohomology ring of the Grassmannian.
The product or {\it intersection} of  any two Schubert classes   
$\sigma_{a_1, a_2}$ and  $\sigma_{b_1, b_2}$
has the form
$$
\sigma_{a_1, a_2}\,\cdot\,\sigma_{b_1, b_2}\,=\,
\sum_{c_1+c_2=a_1+a_2+b_1+b_2} C(a_1,a_2;b_1,b_2; c_1,c_2) 
\sigma_{c_1, c_2}\,,
$$
where $ C(a_1,a_2;b_1,b_2; c_1,c_2)$ are nonnegative integers
called {\it the Littlewood--Richardson coefficients}.

\medskip
If the sum of the codimensions of classes equals $\dim G_2({d+1})=2d-2$,
then their intersection is an integer (identifying the generator of the top
cohomology group $\sigma_{d-1, d-1}\in H^{4d-4}(G_2({d+1}))$ with $1\in\Z$) 
called {\it the intersection number}.

\medskip
When $(a_1, a_2)=(q,0)$, $0\leq q\leq d-1$, 
the  Schubert varieties $\Omega_{q, 0}$  
are called {\it special} and the corresponding cohomology classes 
$\sigma_q\,=\,\sigma_{q, 0}$
are called {\it special Schubert classes}.

\medskip
The  Littlewood--Richardson coefficients appear as well in the decomposition
of the tensor product of two irreducible finite-dimensional $sl_2$-modules into 
the direct sum of  irreducible $sl_2$-modules. This leads to the following 
claim connecting Schubert calculus and representation theory.  

\medskip
Denote $L_q$ the irreducible $sl_2$-module with highest weight $q$.

 \begin{thm}\label{Pr} {\rm \cite{F}}\ \ 
 Let $q_1,\dots, q_{n+1}$ be integers such that   $0\leq q_i\leq  d-1$  
for all $1\leq i\leq n+1$ and  $q_1+\dots +q_{n+1}=2d-2$. The
intersection number of the  special Schubert classes,
$ \sigma_{q_1}\cdot {\rm\ ...\ } \cdot \sigma_{q_{n+1}}\,$,  coincides with  
the multiplicity of  the trivial $sl_2$-module  $L_0$ in
the tensor product  $L_{q_1}\otimes\dots\otimes L_{q_{n+1}}$.  \hfill $\triangleleft$
 \end{thm}

Theorems \ref{tI} and \ref{Pr} imply the following explicit formula.

\begin{thm}  Let $q_1,\dots, q_{n+1}$ be integers such that   $0\leq q_i\leq  d-1$
for all $1\leq i\leq n+1$ and  $q_1+\dots +q_{n+1}=2d-2$. Then 
$$
\sum_{l=1}^n (-1)^{n-l}
\sum_{1\leq i_1<\dots <i_l\leq n}
{q_{i_1}+\dots +q_{i_l}+l-d-1\choose n-2}\,
$$
is the intersection number $\sigma_{q_1}\cdot {\rm\ ...\ } \cdot \sigma_{q_{n+1}}$.
 \hfill $\triangleleft$
\end{thm}
We did not find this formula in the literature on the Schubert calculus.

\subsection{Planes with a given Wronskian  and  Schubert classes}\label{s42}

Applying the Schubert calculus to the Wronski map, 
we arrive at the following result.

\begin{thm}\label{tWr}
Let $m_1\,,\dots, m_n\,$ be positive integers,   $|M|=m_1+\dots+m_n$. 
For generic $z\in Z$  and for any integer $k$ 
such that $1\leq k<|M|+1-k$, the preimage of  the polynomial  {\rm (\ref{W})}  
under the Wronski map consists of at most
 $$
\sigma_{m_1}\,\cdot  {\rm\ ...\ } \cdot\, \sigma_{m_n}\,\cdot \,\sigma_{|M|-2k}
$$
nondegenerate planes of order $k$ and of degree $<|M|+1-k$.
\end{thm}

\noindent{\bf Proof:}\ \  Any plane of  order $k$ with the Wronskian of degree $|M|$
 lies in  $G_2( \Poly_{|M|+1-k})$, as Lemma \ref{l1} shows.
Therefore it is enough to consider the Wronski map on   $G_2( \Poly_{|M|+1-k})$.

\medskip
For  any  $z_j$, define the  flag  $\F_{z_j}$ in $ \Poly_{|M|+1-k}\,$, 
$$
 \F_0(z_j) \subset \F_1(z_j)\subset \dots \F_{|M|+1-k}(z_j),=\, \Poly_{|M|+1-k}\,,
$$
where  $\F_i(z_j)$  consists of the polynomials 
$ P(x)\in\Poly_{|M|+1-k}\,$ of the form
$$
 P(x)= a_i(x-z_j)^{|M|+1-k-i}+\dots +a_0(x-z_j)^{|M|+1-k}\,.
 $$
We have  $\dim\F_{i}(z_j)=i+1$.
Lemma \ref{l4} implies that the nondegenerate planes  with a Wronskian
having  at $z_j$ a root of multiplicity $m_j$ lie in the special
Schubert variety   
$$
\Omega_{m_j,  0}(\F_{z_j})\subset G_2(\Poly_{|M|+1-k})\,.
$$

\medskip
The maximal possible degree of  the Wronskian  $W_V(x)$ for 
$V\in\Poly_{|M|+1-k}$ is clearly $2|M|-2k$.
Denote $m_{\infty}\,=\,|M|-2k$, the difference between $2|M|-2k$ and $|M|$.
If $m_{\infty}$ is positive, it is
{\it the multiplicity of $W(x)$ at infinity}. 

\medskip
The nondegenerate planes  with a Wronskian
having  given multiplicity $m_{\infty}$ at infinity lie in   
the special Schubert variety   
$\Omega_{m_{\infty}, 0}(\F_\infty)\,$, where $\F_\infty$ is the flag
$$
\Poly_0\subset\Poly_1\subset\dots\subset \Poly_{|M|+1-k}\,.
$$

\medskip
We conclude that the nondegenerate planes of order $k$ which have the Wronskian
 (\ref{W})  lie in the intersection of special 
Schubert varieties
$$
\Omega_{m_1, 0}(\F_{z_1})\,\cap \Omega_{m_2, 0}(\F_{z_2})
\cap \dots \cap \Omega_{m_n, 0}(\F_{z_n})
\, \cap \Omega_{m_{\infty}, 0}(\F_{\infty})\,.
$$
\medskip
The dimension of  $G_2(\Poly_{|M|+1-k})$  is exactly
$m_1+\dots+m_n+m_{\infty}\,$, therefore  this intersection 
consists of a finite number of planes, and
the intersection number of the special Schubert classes 
$$
\sigma_{m_1}\,\cdot {\rm\ ...\ }\cdot \sigma_{m_n}\,\cdot\, \sigma_{|M|-2k}\,
$$
provides an upper bound.   \hfill $\triangleleft$    

\medskip
Theorems \ref{Pr} and  \ref{tWr}, together with Corollary \ref{c2},  
imply the Proposition  of  Sec.~\ref{s27}. 

\begin{cor}\label{last}
For generic $z\in Z$, all planes in the preimage under the Wronski map
of the polynomial $W(x)$ given by {\rm (\ref{W})}  are nondegenerate.
\hfill $\triangleleft$
\end{cor}

\section{The simplicity of the spectrum of the Gaudin hamiltonians}\label{S5}
 \subsection{Fuchsian differential equations with only polynomial
solutions and the Wronski map}\label{s51}
Consider a second order Fuchsian differential equation with 
regular singular points at $z_1,\dots,z_n$, $n\geq 2$, and at 
infinity.  If   the exponents at  $z_j\,$ are $0$ and $m_j+1$,
$1\leq j\leq n$, then this equation has  the form
\begin{eqnarray}\label{FE}
F(x)u''(x)+G(x)u'(x)+H(x)u(x)=0\,,\\
F(x) = \prod_{j=1}^n (x-z_j)\,,\quad
\frac{G(x)}{F(x)} = \sum_{j=1}^n\frac{-m_j}{x-z_j}\,,\nonumber
\end{eqnarray}
where $H(x)$ is a polynomial of degree not greater than $n-2$.
On Fuchsian equations see Ch.~6 of \cite{R}.

\begin{thm}
Any nondegenerate plane $V\in G_2(\Poly_d)$  with the Wronskian {\rm (\ref{W})}
is the solution space of  such an equation.
\end{thm}

\noindent{\bf Proof:}\ \ 
Let  $f(x)$ and  $g(x)$ be non-zero monic polynomials  in $V$ such that $\deg f<\deg g$.
As we have seen in course of the proof of Lemma~\ref{l2}, the plane $V$ is the
solution space of the equation
$$
W(x)u''(x)-W'(x)u'(x)+h(x)u(x)=0\,,
$$
where 
$$
h(x)\, =\, \frac{-W(x)f''(x)+W'(x)f'(x)}{f(x)}\,
$$
is a polynomial  proportional to the Wronskian of
$f'(x)$ and $g'(x)$. One can easily check that
$$
\frac{W'(x)}{W(x)}\,=\,\sum_{j=1}^n\frac{m_j}{x-z_j}\,.
$$
Moreover, if $z_j$ is a root of $W(x)$ of multiplicity $m_j>1$,
then all coefficients of the equation have $(x-z_j)^{m_j-1}$ as
a common factor, and the equation can be reduced to the required  form (\ref{FE}).  
\hfill $\triangleleft$

\medskip
In \cite{Fr}, a correspondence between  eigenvalues of the Gaudin hamiltonians 
and the $sl_2$-opers which determine a trivial monodromy representation
$$
\pi_1(P^1\setminus\{z_1,\dots,z_n,\infty\}) \to PGL_2
$$ 
was established. Notice that  
any nondegenerate plane with  a given  Wronskian defines such an oper,
and  any oper with trivial monodromy  defines a plane in $ G_2(\Poly_d)$
for some $d$.

\subsection{Bethe vectors are differ by their eigenvalues}\label{s52} 
Let $t^{(1)}$ and $t^{(2)}$ be two solutions to the Bethe equations (\ref{cps}) 
associated to fixed $z^0\in Z$ and $\Sing_k$. Denote $v(t^{(i)},z^0)$ the corresponding 
Bethe vectors given by   (\ref{Bv}) and  $\mu^{(i)}=(\mu_1(t^{(i)}),\dots, \mu_n(t^{(i)}))$ 
their eigenvalues given by (\ref{cvs}), $i=1,2$.

\begin{thm}
If  $\mu^{(1)}=\mu^{(2)}$, then  $v(t^{(1)},z^0)= v(t^{(2)},z^0)$.
\end{thm}

\noindent{\bf Proof:}\ \ 
Let $f_i(x)=(x-t^{(i)}_1)\dots (x-t^{(i)}_k)$ be the corresponding polynomials,  $i=1,2$.
Any of these polynomials defines  a differential equation of the form  (\ref{FE}).
A classical fact of the theory of Fuchsian equations is that if  an equation of the 
form  (\ref{FE}) with positive integer $m_1,\dots,m_n$ has a polynomial solution without 
multiple roots, then all  solutions to this equation are polynomials, see Ch.~6 of \cite{R} 
or  Sec.~3.1 of \cite{SV}. We have 
$$
H_i(x)=-\frac {F(x)f_i''(x)+G(x)f_i'(x)}{f_i(x)}\,,\quad i=1,2.
$$
Polynomials $H_1(x)$ and $H_2(x)$ have  degree at most $n-2$. 
We will show that these polynomials in fact coincide. Indeed,
the substitution of $z_j$  into $H_i$ gives
$$
H_i(z_j)=m_j \cdot\, \frac{f_i'(z_j)}{f_i(z_j)}\,\prod_{l\neq=j}^n(z_j-z_l).
$$
An easy calculation shows that
$$
\frac{f_i'(z_j)}{f_i(z_j)}=\sum_{l=1}^k\frac1{z_j-t^{(i)}_l},
$$
and the condition $\mu^{(1)}=\mu^{(2)}$ together with (\ref{cvs}) imply
$$
\frac{f_1'(z_j)}{f_1(z_j)}=\frac{f_2'(z_j)}{f_2(z_j)}\,,\quad j=1,\dots,n.
$$
We conclude that  polynomials $H_i(x)$  coincide at $n$ points, therefore 
$H_1(x)=H_2(x)$.  Thus polynomials $f_1(x)$ and $f_2(x)$
are monic polynomials of  the minimal degree in the same solution plane. Hence
$f_1(x)=f_2(x)$,  the solutions $t^{(1)}$ and $t^{(2)}$ of the Bethe equations
lie in the same orbit, and  $v(t^{(1)},z^0)= v(t^{(2)},z^0)$.
\hfill $\triangleleft$

\begin{cor}
For generic $z^0\in Z$, the $sl_2$ Gaudin hamiltonians {\rm (\ref{H})} have  
a simple spectrum. \hfill $\triangleleft$
\end{cor}

\bigskip

\

\end{document}